\documentclass[12pt]{amsart}
\usepackage{amsfonts, amssymb,amsmath, amsthm, amsxtra, latexsym, mathrsfs, amscd}
\usepackage{amsmath, amssymb, amsthm, times, float}
\usepackage{mathtools, amssymb, amsthm, times, float}
\usepackage{graphics}
\usepackage{MnSymbol}

\usepackage{amssymb,hyperref}
\usepackage{backref}
\usepackage[latin1]{inputenc} %

\newtheorem{theo+}{Theorem}[section]
\newtheorem{prop+}[theo+]{Proposition}
\newtheorem{coro+}[theo+]{Corollary}
\newtheorem{lemm+} [theo+]{Lemma}
\newtheorem{deep+}  [theo+]  {Deep Result}
\newtheorem{fact+}  [theo+]  {Fact}
\theoremstyle{definition}
\newtheorem{exam+}  [theo+]  {Example}
\newtheorem{rema+}  [theo+]  {Remark}
\newtheorem{defi+}  [theo+]  {Definition}
\newtheorem{xca+}[theo+]{Exercise}

\numberwithin{equation}{section}

\usepackage{color}

\newcommand\beq{\begin{equation}\label}
\newcommand\eeq{\end{equation}}

\renewcommand\){\Big)}

\renewcommand\a[1]{{\acute{#1}}}

\def\draft{\centerline{(Draft {\the \day}/{\the\month} \the \year.)}}
\bigskip

\def\refn#1.#2{\expandafter\def\csname#1\endcsname{[#2]}}
\def\refnr#1.{\csname#1\endcsname}

\def\fg{\mathfrak g}
\def\fgl{\mathfrak {gl}}
\def\fk{\mathfrak k}

\def\fp{\mathfrak p}

\def\fsl{\mathfrak {sl}}

\def\fu{\mathfrak u}

\def\fsp{\mathfrak{sp}}

\def\fsu{\mathfrak{su}}

\def\a{\alpha}

\def\Claminv2{|C(\Lambda)|^{-2}}

\def\de{d\varepsilon}

\def\Aa2D{A^{\a,2}(D)}
\def\bAa2D{\overline{A^{\a,2}(D)}}
\def\Ab2D{A^{\beta,2}(D)}
\def\bAb2D{\overline{A^{\beta,2}(D)}}

\def\Norm#1_#2{\Vert#1\Vert_{#2}}

\def\phipl12{\phi_{p_{l_1}, p_{l_2}}}
\def\phip01{\phi_{p_{0}, p_{0}}}

\def\a{\alpha}

\def\Claminv2{|C(\Lambda)|^{-2}}

\def\sig{\sigma}

\def\det{\operatorname{det}}
\def\diag{\operatorname{diag}}

\def\diag{\operatorname{diag}}

\def\exp{\operatorname{exp}}

\def\tr{\operatorname{tr}}

\def\Ker{\operatorname{Ker}}

\def\de{d\varepsilon}

\def\Aa2D{A^{\a,2}(D)}
\def\bAa2D{\overline{A^{\a,2}(D)}}
\def\Ab2D{A^{\beta,2}(D)}
\def\bAb2D{\overline{A^{\beta,2}(D)}}

\def\phipl12{\phi_{p_{l_1}, p_{l_2}}}
\def\phip01{\phi_{p_{0}, p_{0}}}

\def\alg/{algebra}

\def\Alg/{Algebra}

\def\alt/{alternative} 
\def\anal/{analytic}
\def\analfunc/{\anal/\ \func/}
\def\Ans/{\it Answer. \normal}
\def\ass/{associative}
\def\nass/{non-\ass/}
\def\autom/{automorphism}
\def\homom/{homomorphism}
\def\isom/{isomorphism}
\def\bdd/{bounded}
\def\Bdd/{Bounded}
\def\bddsymdom/{bounded \sym/ \dom/}
\def\Cartdom/{Cartan \dom/}
\def\bdry/{boundary}
\def\bsd/{\bdd/ \symdom/}
\def\bv/{boundary value}
\def\cf/{{\it cf}\.}
\def\Cf/{{\it Cf}\.}
\def\charr/{character}
\def\coeff/{coefficient}
\def\comm/{commutative}
\def\cpct/{compact}
\def\compl/{complex}
\def\comp/{complex}
\def\Comp/{Complex}
\def\conf/{conformal}
\def\conj/{conjugate}
\def\conn/{connect}
\def\cont/{continuous}
\def\conv/{converge} 
\def\convc/{convergence}
\def\convt/{convergent}
\def\convx/{convex}
\def\coord/{coordinate}
\def\lcoord/{local coordinate}
\def\Corr/{Corresponding}
\def\corr/{corresponding}
\def\corrd/{correspond}
\def\cov/{covariant}
\def\decomp/{decomposition}
\def\deco/{decompose}
\def\diff/{different} 
\def\Diff/{Different} 
\def\dimn/{dimension} 
\def\distr/{distribution} 
\def\div/{diverge} 
\def\dom/{domain}
\def\eg/{\hbox{\it e.g}\.}
\def\eigenf/{eigen\-\func/}
\def\eigensp/{eigen\-space}
\def\eigenv/{eigen\-value}
\def\eq/{equation}
\def\equa/{equation}
\def\de/{\diff/ial \equa/}
\def\do/{\diff/ial operator}
\def\ode/{ordinary \de/}
\def\pde/{partial \de/}
\def\pdo/{partial \diff/ial operator}
\def\psdo/{pseudo \diff/ial operator}
\def\fin/{finite}
\def\Ex/{\it Example.\ \normal}
\def\Exnr#1/{\it Example #1.\ \normal}
\def\foll/{follow}
\def\follg/{following}
\def\Follg/{Following}
\def\func/{function}
\def\Func/{Function}
\def\Fonc/{Fonc\-tion}
\def\fonc/{fonc\-tion}
\def\Funk/{Funk\-tion}
\def\funk/{Funk\-tion}
\def\gen/{general}
\def\har/{harmonic}
\def\Hint/{\it Hint. \normal}
\def\hist/{historic}
\def\histcl/{historical}
\def\hol/{holo\-morphic}
\def\homog/{ho\-mo\-ge\-ne\-ous}
\def\hyp/{hyper\-bolic}
\def\hyperg/{hyper\-geometric}
\def\ie/{\hbox{\it i.e.}}
\def\iff/{if and only if}
\def\ineq/{inequality}
\def\infra/{{\it inf\-ra}}
\def\ultra/{{\it ult\-ra}}
\def\Inpart/{In particular}
\def\inpart/{in particular}
\def\instof/{instead of}
\def\interps/{interpolation space}
\def\interp/{interpolation}
\def\Interp/{Interpolation}
\def\interpr/{Interpretation}
\def\Intr/{Introduction}
\def\intv/{interval}
\def\inv/{invariant}
\def\invc/{invariance}
\def\Iowords/{In other words}
\def\iowords/{in other words}
\def\ipr/{inner product}
\def\irred/{irreducible}
\def\lb/{line bundle}
\def\lin/{linear}
\def\lhs/{left hand side}
\def\rhs/{right hand side}
\def\loc/{local}
\def\math/{mathematic}
\def\mathcn/{\math/ian}
\def\manif/{manifold}
\def\meas/{measure}
\def\measl/{measurable}
\def\mero/{mero\-morphic}
\def\mon/{monomial}
\def\monog/{monogenic}
\def\mult/{multiple}
\def\multy/{multiply}
\def\multn/{multiplication}
\def\nas/{necessary and sufficient}
\def\nbd/{neighborhood}
\def\neg/{negative}
\def\nondeg/{nondegenerate}
\def\Oohand/{On the other hand}
\def\oohand/{on the other hand}
\def\Oonhand/{On the one hand}
\def\oonhand/{on the one hand}
\def\oper/{operator}
\def\orth/{ortho\-gonal}
\def\orthon/{ortho\-normal}
\def\otoh/{on the other hand}
\def\quat/{quaternion}
\def\pp/{\hbox{a. e.}}
\def\psh/{plurisubharmonic}
\def\pol/{polynomial}
\def\pot/{potential}
\def\pos/{positive}
\def\princ/{principle}
\def\prob/{probability}
\def\proj/{projective}
\def\projn/{projection}
\def\Proof/{\it Proof:\normal}
\def\Rem/{\it Remark\normal}
\def\Remnr#1/{\it Remark\ \normal #1. }
\def\rep/{representation}
\def\reps/{representations}
\def\meta/{metaplectic representation}
\def\repr/{reproducing}
\def\reprker/{reproducing kernel}
\def\resp/{respective} 
\def\resply/{respectively}
\def\restr/{restriction}
\def\sa/{self-adjoint}
\def\st/{such that}

\def\sol/{solution}
\def\ru/{space}
\def\sph/{spherical}
\def\ssp/{sub\ru/}
\def\sym/{symmetric}
\def\Sym/{Symmetric}
\def\symb/{symbol}
\def\symbc/{symbolic}
\def\symdom/{\sym/ domain}
\def\symp/{symplectic}
\def\Theor#1/{\fet Theorem #1.\ \normal}
\def\Lem#1/{\fet Lemma #1.\ \normal}
\def\Lemma/{\fet Lemma.\ \normal}
\def\topl/{topology}
\def\topll/{topological}
\def\transf/{transform}
\def\transl/{translation}
\def\transfn/{transformation}
\def\transv/{transvectant}
\def\trig/{trigonometric}
\def\tril/{trilinear}
\def\trilf/{trilinear form}
\def\uhp/{upper halfplane}
\def\uhs/{upper halfspace}
\def\vb/{vector bundle}
\def\vf/{vector field}
\def\vsp/{vector space}
\def\wrt/{with respect to}
\def\Wlog/{Without loss of generality}
\def\a{\alpha}

\def\sig{\sigma}

\def\Ab/{Abel}
\def\Ban/{Banach}
\def\Bansp/{\Ban/ space}
\def\Belt/{Bel\-tra\-mi}
\def\Berg/{Berg\-man}
\def\Bern/{Ber\-nou\-lli}
\def\Berz/{Berezin}
\def\Bess/{Bessel}
\def\Cart/{Car\-tan}
\def\Cay/{Cay\-ley}
\def\CG/{Clebsch-Gordan}
\def\Cl/{Clifford}
\def\CR/{Cauchy-Rie\-mann}
\def\Dir/{Dirichlet}
\def\Eucl/{Euclide}
\def\Eucln/{Euclidean}
\def\F/{Fourier}
\def\Hank/{Hankel}
\def\Hankf/{\Hank/ form}
\def\Herm/{Hermite}
\def\Hilb/{Hilbert}
\def\Hilbs/{Hilbert space}
\def\Hilbsp/{Hilbert space}
\def\HS/{Hilbert-Schmidt}
\def\Lag/{La\-grange}
\def\Lap/{La\-place}
\def\LapBelt/{\Lap/-\Belt/}
\def\Leb/{Lebesgue}
\def\Marc/{Mar\-cin\-kie\-wicz}
\def\Moeb/{Moebius}
\def\Moebt/{Moebius transformation}
\def\Moebtransfn/{Moebius transformation}
\def\Pla/{Plan\-che\-rel}
\def\Poin/{Poin\-car\'e}
\def\Riem/{Rie\-mann}
\def\Riemn/{\Riem/ian}
\def\psRiemn/{pseudo-\Riem/ian}
\def\Riems/{Rie\-mann surface}
\def\Schroe/{Schr\"odinger}
\def\Weier/{Weier\-strass}
\def\anal/{analytic}
\def\bsd/{bounded symmetric domain  }
\def\bdd/{bounded}
\def\calc/{calculation}\def\conj{conjugate}
\def\calci/{calculating}\def\eg{e.g.}
\def\conj/{conjugate}
\def\deco/{decomposition}
\def\eg/{e.g.}
\def\fct/{function}
\def\gp/{group}
\def\hw/{highest weight}
\def\hwv/{highest weight vector}
\def\hwvs/{highest weight vectors}
\def\lw/{lowest weight}
\def\lwv/{lowest weight vector}
\def\lwvs/{lowest weight vectors}
\def\hds/{holomorphic discrete series}
\def\iff/{if and only if}
\def\inv/{invariant}
\def\irrde/{irreducible decomposition}
\def\meas/{measure}
\def\transf/{transform}
\def\rep/{representation}
\def\resp/{respectively}
\def\inters/{intertwines}
\def\interg/{intertwining}
\def\meta/{metaplectic representation}
\def\qu/{quaternion}
\def\rep/{representation}
\def\symdom/{ symmetric domain}
\def\st/{such that}
\def\shd/{subhead}
\def\transf/{transform}
\def\wrt/{with respect to}

\def\Norm#1#2#3{\Vert#1\Vert^{#3}_{{#2}}}

\def\tr{\operatorname{tr}}

\begin{document}

\title[
Branching of metaplectic representation 
of $Sp(2, \mathbb R)$ under its principal 
$SL(2, \mathbb R)$-subgroup
]
{
Branching of metaplectic representation  
of $Sp(2, \mathbb R)$ under its principal  
$SL(2, \mathbb R)$-subgroup 
}
\author{GenKai Zhang}
\address{Mathematical Sciences, Chalmers University of Technology and
Mathematical Sciences, G\"oteborg University, SE-412 96 G\"oteborg, Sweden}
\email{genkai@chalmers.se}

\thanks{ Research supported partially
 by the Swedish
Science Council (VR). 
}
\begin{abstract}
  We study the  branching problem
of  the metaplectic representation of $Sp(2, \mathbb R)$
under its  principle subgroup $SL(2, \mathbb R)$.
We find  the complete decomposition.
\end{abstract}

\maketitle

\baselineskip 1.35pc

\section{Introduction}

The problem of finding the branching
rule of a unitary representation $(\Lambda, G)$
of a semisimple Lie group $G$
under a subgroup $H$ is
of fundamental interests in representation theory.
In the present paper we shall study
this problem for  the metaplectic
representation $\Lambda$
of the symplectic group $G=Sp(2, \mathbb R)$
under its principal subgroup $SL(2, \mathbb R)$.

The metaplectic 
representation for $G=Sp(2, \mathbb R)$
can be realized on the Fock space $\mathcal F(\mathbb C^2)$
and the Lie algebra $\fg$ acts as 
differential operators with quadratic
polynomial coefficients.
We use the classical  approach
by finding the spectra of the Casimir operator $C$.
To find the continuous spectrum
  we compute the action of the Casimir element  
  on some orthogonal basis of  weight vectors
  of $SO(2)\subset SL(2, \mathbb R)$  of a fixed weight.
We prove that Casimir element is 
 unitary equivalent
   to the coordinate multiplication operator on
   certain hypergeometric orthogonal polynomials,
   more precisely the continuous dual Hahn polynomials, and it
   determines
   the continuous spectrum of $C$. 
   
   To find the discrete components we solve the equation
   of highest or lowest weight vectors for $\fsl(2, \mathbb C)$-actions in the Fock space.
   We prove    that lowest weight representations
   of $SL(2, \mathbb R)$
   will not appear in $(\Lambda, G)$ and we find all highest weight
   representations.

For the even part 
   of the metaplectic representation 
   this decomposition can also be treated 
   using the Berezin transform \cite{gz-bere-rbsd}; see 
   Remark \ref{rema-be}. 
   The connection between 
   hypergeometric  orthogonal polynomials and branching 
   problems has been observed  in the earlier works \cite{gz-brc-jfa,
     gz-rep-thy},
   it was proved that density for orthogonality is precisely the product 
   of Harish-Chandra $c$-function and the symbol of the Berezin transform.

     We remark that if  we transfer our results 
   in the Schr\"o{}dinger  $L^2(\mathbb R^2)$-realization, we have 
   essentially found the spectrum of 
some   partial differential operators 
of degree four involving squares of harmonic oscillator 
of the form  $\alpha (A\partial, \partial) + \beta (Bx, x)$
   with non-positive definite symmetric matrices $A$ and $B$, and 
   this might be of independent interests. 

   Finally we mention very briefly
   our motivation, related known results
   and some perspective questions.
   Firstly the exact appearance of $SL(2, \mathbb R)\subset G=Sp(2, \mathbb R)$
   and the corresponding branching problem
   can be used to study the induced representation
   of the group $G_2$ from its Heisenberg
   parabolic subgroup \cite{Frahm};
   in \cite{FWZ} we shall study the non-compact
   relization of induced
   represenations for Hermitian Lie groups
   using the related branching of metaplectic
   representation, the compact picture being studied in
   \cite{Z-21}. Next this kind
   of branching problem can be formulated for any  real split simple Lie group $G$. Indeed
up to conjugation  
there is a unique principal subalgebra $\fsl(2, \mathbb R)$ (and
it might have more than one principal $\fsl(2, \mathbb R)$-subgroups); see \cite{Kostant}.
   For  general  $G=Sp(n, \mathbb R)$
   we can use covariant differentiations to produce 
some discrete components, but the results 
are not complete as for 
$G=Sp(2, \mathbb R)$ and will not be present here; 
see  e.g. \cite{DGV,
Ko,  Orsted-Vargas} 
for related study and references.
It is also tempting but
challenging to study the branching 
problems for the minimal representations of $G$
under $SL(2, \mathbb R)$, even for groups of lower
rank such as $SL(3, \mathbb R)$ and $G_2$.
The principal $SL(2, \mathbb R)$-subgroups
in split real Lie
groups corresponde to certain nilpotent orbits. On the other
hand general nilpotent orbits of real Lie groups
are closely related  $SL(2, \mathbb R)$-subgroups.
One might ask if there is a general theory
about the branching under  $SL(2, \mathbb R)$
of representations of $G$
related to the corresponding nilpotent orbits. Next
the branching of metaplectic representations
under  dual pairs, namely the dual correspondence,
has been under intensive study and has found
many applications; for the dual pair $SL(2, \mathbb R)\times O(p, q)$
in $Sp(p+q, \mathbb R)$ there
is some detailed study of Plancherel formula
by Howe \cite{Howe}  using the results
of Repka \cite{Repka} on tensor products of
representations of  $SL(2, \mathbb R)$; however our result
here seems can not be deduced from
\cite{Repka}. From a classical analytic
point of view explicit decompositions under $SL(2, \mathbb R)$
provide methods to discover new orthogonal polynomials
and spectral decomposition of higher order self-adjoint operators.
Finally
the appearance of certain
specific  $SL(2, \mathbb R)$
in $Sp(n, \mathbb R)$ and hyperbolic discs
in locally symmetric spaces or Teichm\"uller
space is of importance in higher Teichm\"u{}ller theory \cite{BIW}.

   The paper is organized as follows.  In Section 1 we recall the realization
   of the metaplectic representation of $Sp(n, \mathbb R)$
   on the Fock space.
We find the explicit description 
  of the principal $SL(2, \mathbb R)$-subgroup in $Sp(n, \mathbb R)$
  and its action on the metaplectic representation
  in Section 2. The complete  decomposition
  of the metaplectic representation of
  $Sp(2, \mathbb R)$ is done in Section 3.

{\bf Acknowledgement.}
I would like to thank  Wee Teck Gan for posing
the question studied in this paper  and for
some stimulating discussions. I'm grateful to
the anonymous referee for a careful
reading of an earlier version
of this paper and for several expert suggestions.

\section{Preliminaries 
}

We recall some known facts on metaplectic representations
of $G=Sp(n, \mathbb R)$.

\subsection{The group $G=Sp(n,  \mathbb R)$
}
\label{1.1}
The symplectic group
$G=Sp(n, \mathbb R)$
has usually two different realizations,
 as the group of biholomorphic transformations
of the Siegel upper half space
or  the bounded symmetric domain.
We shall  use the latter one, which is sometimes denoted
by $G=Sp(n, \mathbb R)_c$. See \cite{FK-book, FK, Folland, Loos-bsd}.

Let
$Sp(n, \mathbb C)$ be the group
of complex $2n\times 2n$-matrices preserving the
standard complex symplectic  form
$\Omega(z, w)=z_1w_{n+1}+\cdots + z_n w_{2n} - 
z_{n+1}w_1-\cdots - z_{2n} w_n.$
Let $U(n, n)$ be the indefinite unitary group
on $\mathbb C^{2n}$ preserving the Hermitian
form $\langle z, w\rangle_{n, n}=z_1\overline w_1
+\cdots + z_n\overline w_n-   \cdots- z_{2n}\overline w_{2n}$ of signature
$(n, n)$.
Let  $G=Sp(n, \mathbb R):=U(n, n)\cap Sp(n, \mathbb C).$
Elements in $G$
will be represented as $2\times 2$-block matrices
\begin{equation}
  \label{g-a-b}
g=\begin{bmatrix} a & b\\
  \bar b & \bar a
 \end{bmatrix}, 
\end{equation}
where    $a, b $ are 
complex    $n\times n$-matrices.
Let $K\simeq U(n)$ be the subgroup of
 block diagonal matrices
 $g=\diag(a, \bar a)$, with $a\in U(n)$,
 $U(n)$ being in its standard realization.
 Then $K$ is the maximal subgroup of $G$,
 and to save notation we
 identify  $\diag(a, \bar a)$ with $a$.

  The Cartan decomposition  is
  $\fg=\fk+\fp$ with $\fk=\fu(n)$, $\fp=M_{n}^s(\mathbb C)
  =\{B\in M_n(\mathbb C); B=B^t\}$; more precisely 
   \begin{equation}
     \label{cartan-k}
     \mathfrak k=\left\{diag(A, -A^t), A\in \fu(n)\right\},
   \end{equation}
   \begin{equation}
     \label{cartan-p}
  \mathfrak p=\left\{\xi_B=\begin{bmatrix} 0 &B\\
  \bar{B} &0
\end{bmatrix}; B=B^t\in M_n(\mathbb C)\right\},
   \end{equation}
in the above matrix realization.

\subsection{Metaplectic  representation}

 Let $\mathcal F=\mathcal F(\mathbb C^n)$
 be the Fock space of entire functions $f\in \mathcal O(\mathbb C^n)$
 on
$\mathbb C^n$ such that
$$
\Vert f\Vert^2=\int_{\mathbb C^n}|f(z)|^2 e^{-\pi |z|^2}
dz<\infty,
$$
where $dz$ is the Lebesgue measure.
The monomials $\{z^{\alpha}
=
z_1^{\alpha_1}\cdots z_n^{\alpha_n}\}$
form an orthogonal basis with norm square
\begin{equation}
  \label{monom-norm}
  \Vert
  z^{\alpha}\Vert^2
  =
  \frac 1{\pi^{\alpha_1 +\cdots +
      \alpha_n}}
  (\alpha_1!)\cdots (\alpha_n!),
\end{equation}
and  $\mathcal F$ has  reproducing kernel  $e^{\pi (z, \bar w)}$, namely
\begin{equation}
  \label{rep-k}
  f(z)=\langle f, e^{\pi (\cdot, \bar z)}
  \rangle=\int f(w)  e^{\pi (z, \bar w)} e^{-\pi |w|^2} dw, \quad f\in \mathcal F, \,
z\in \mathbb C^n,  
\end{equation}
where  $(z, w)=\sum z_j w_j$ is the
$\mathbb C$-bilinear form.

The metaplectic representation $(\Lambda, \mathcal F, Mp(n, \mathbb R))$
is
 a unitary 
representation of the 
double cover  of $G=Sp(n, \mathbb R)$
on the Fock space $\mathcal F$,
and is explicitly  given by
 \cite[Theorem 4.37]{Folland},
$$
\Lambda(g) f(z)=\det^{-\frac 12}(a) \int
\exp
\left\{
    \frac \pi 2 \left[(\bar b a^{-1} z, z)
+ 2
(a^{-1} z, \bar w) - (a^{-1}b\bar w, \bar w)
\right]\right\} f(w) e^{-\pi |w|^2}dw, 
$$
for $g$ as in
(\ref{g-a-b}).
In particular the group $U(n)$ acts as
$$
\Lambda (a)f(z)= \det^{-\frac 12}(a)
f(a^{-1}z), \quad 
\quad a\in U(n).
$$


The space $\mathcal F=\mathcal F_+\oplus \mathcal F_-
$ is  sum of the subspaces of even and odd functions,
with $\Lambda=\Lambda^+\oplus \Lambda^-$
a sum of two irreducible representations.

\begin{lemm+}
\label{Lie-alg-mu}
  The Lie algebra action of $sp(n, \mathbb C)$
 on $\mathcal F$ is given as follows
  $$
  \Lambda (X)f =
  -\frac{1}{2\pi}
  (B\partial, \partial)f, \quad
  \Lambda (Y)f =\frac \pi 2(\overline B z, z)f, 
  $$
$$
  \Lambda (Z)f (z)=
  -\frac{\tr D}{2}f(z)
 -(\partial_{Dz}f)(z),
  $$
  for
  $$
  X=X_B=
  \begin{bmatrix} 0 & B\\
    0 & 0
    \end{bmatrix}, \quad  Y=Y_{\overline  B}=
  \begin{bmatrix} 0 & 0\\
    \overline   B & 0
  \end{bmatrix}, \quad B=B^t, 
  $$
  $$
Z=Z_D=  \begin{bmatrix} D & 0\\
     0 & -D^t
    \end{bmatrix}, \quad D\in  \fgl(n).
$$
  \end{lemm+}
  \begin{proof} Let $\xi_B$
    be as in          (\ref{cartan-p}) and perform
    the differentiation $\frac{d}{dt}$ of the action $\exp(t \xi_B)\in G$
    at $t=0$,
    $$
    \Lambda(\xi_B) f(z)
    =\int \frac \pi 2 \left[
        (\overline B z, z) - (B\bar w, \bar w)\right]
    e^{\pi ( z, \bar w)}
    f(w) e^{-\pi |w|^2}dw.
  $$
  The  integration of the first term
  is
  $$
  \frac \pi 2 (\overline B z, z) f(z)  
  $$
  by the reproducing kernel
  formula (\ref{rep-k}).
  Differentiating the formula
  (\ref{rep-k})
  by $(B\partial, \partial)$
  we find
  $$
  (B\partial, \partial) f(z)
  =\pi^2 \int (B\bar w, \bar w)
      e^{\pi ( z, \bar w)}  f(w) e^{-\pi |w|^2}dw.
  $$
  Thus the second integration above is 
  $$
  -  \frac {\pi} 2 \frac {1}{\pi^2}  (B\partial, \partial) f(z)
  =  -   \frac {1}{2 \pi}  (B\partial, \partial) f(z).
  $$
    Taking the complex linear and conjugate linear
    parts in $B$ we get the first two formulas.
   The last formula is a straightforward computation.
\end{proof}

\subsection{Unitary representations of $SU(1, 1)$ and their realizations}
\label{1.3}
Let  $G_0:=SU(1, 1)$ be the group
of  $2\times 2$-matrices preserving
the Hermitian form $\langle z, w\rangle_{1, 1}=z_1\bar w_1 -z_2\bar w_2$
with determinant $1$.
We fix the maximal subgroup
\begin{equation}
  \label{eq:U(1)}
L=U(1)=\left\{u_\theta=\begin{bmatrix}
    e^{i\theta} &0\\
0    &
    e^{-i\theta} 
    \end{bmatrix}
\right\}\subset G_0.
\end{equation}
We fix also  the Lie algebra elements 
\begin{equation}
  \label{h-e-pm}
h=\begin{bmatrix} 1& 0\\ 0& -1 
\end{bmatrix},\quad
e^+=\begin{bmatrix} 0& 1\\ 0& 0 
\end{bmatrix},\quad
e^-=\begin{bmatrix} 0& 0\\ 1& 0 
\end{bmatrix}, \quad e=e^{+} +e^{-}
\end{equation}
forming the standard basis of $sl(2, \mathbb C)$. 
Let
$$C=e^+e^-+e^-e^++\frac 12 h^2 
    =2e^+e^- -h +\frac 12 h^2 
    $$
    be the Casimir element.
    
Recall \cite[Theorem 1.3]{HT} that any  unitary irreducible
  representation  $\sig$
  of $SU(1, 1)$ 
  is one of the following 
  \begin{enumerate}
  \item Spherical principal or complementary series
    $\sig_{i\lambda, +}$
    with    all weights of $h$  being even integers;
    the eigenvalue of Casimir $C$
    is $-\frac 12 -\lambda^2$
    for $\lambda \in \mathbb R^+$ or $\lambda \in (0, \frac 12)$.

  \item Non-spherical principal  series
    $\sig_{i\lambda, -}$;
      for $\lambda \in \mathbb R^+$
      with all weights of $h$ are odd integers;
          the eigenvalue of Casimir $C$
          is $-\frac 12 -\lambda^2$
              for $\lambda \in \mathbb R^+$.

            \item Highest weight representations $\sig_{-\nu}$
              with negative integral weights $-\nu$, $\nu \ge 1$;
                        the eigenvalue of Casimir $C$ is $\frac 12 \nu^2 -\nu$.

      \item Lowest weight representations $\sig_{\nu}$ with positive
        integral weights 
        $\nu$, $\nu \ge 1$;
                                the eigenvalue of Casimir $C$ is $\frac 12 \nu^2 -\nu$.
                              \end{enumerate}

                              The non-integral highest or lowest
                              weights $\nu$ corresponde
                              representations of the
                              universal covering of $SU(1, 1)$, which
                              will not concern us here.

\subsection{Continuous 
  dual Hahn polynomials 
  and spectra of multiplication operators 
}

We recall the orthogonality relation
of continuous dual Hahn polynomials. They
are special cases 
of Wilson hypergeometric polynomials
\cite{Koe-Swa},
and will be used in Section 3 below
to find the spectra of Casimir element.

\begin{prop+} (\cite[pp.~196-199]{Koe-Swa})
  Let $a\ge 0, b>0, c>0$, 
$\mu=\mu_{a, b, c}$ be  the measure 
  $$d\mu_{a, b, c}(x) 
  =\frac 1{2\pi}\frac 1{\Gamma(a+b) 
    \Gamma(a+c)    \Gamma(b+c) 
  }
\left|  \frac{\Gamma(a+ix) 
    \Gamma(b+ix) 
    \Gamma(c+ix)    
  }
  {\Gamma(2ix)} \right|^2 
  dx 
  $$
  on $\mathbb R^+=[0, \infty)$ and $L^2(\mathbb R^+, \mu)$
  the corresponding $L^2$-space. 
  Let $$\omega_m(x^2)=\omega_{m; a, b, c}(x^2)=
  {}_3F_2 (-m, a+ix, a-ix; a+b, a+c, 1), 
  \quad m=0, 1, \dots, 
  $$
be the continuous dual Hahn polynomials. Then 
$\{\omega_m(x^2)\}$ form an orthogonal basis of 
$L^2(\mathbb R^+, \mu)$, 
   $$
  \langle  \omega_m, \omega_{l} 
  \rangle 
  =\frac{m!(b+c)_m }{(a+b)_m(a+c)_m} \delta_{m,l}, 
  $$
  where
  $$
  (a)_m = a (a+1)\cdots (a+m-1)
  $$
  is the Pochhammer symbol.
The multiplication operator by $-x^2$ has the following 
$3$-term recursion  on the basis, 
  $$
  -x^2  \omega_{m}(x^2)  =A _m  \omega_{m+1}(x^2) 
  -(A_m +C_m -a^2)  \omega_{m}(x^2) 
  +C_m   \omega_{m-1}(x^2), 
  $$
  where   $$
  A_m =(m+a+b)(m+a+c), \, 
   C_m =m(m+b+c-1). 
  $$
  In particular the spectrum of 
  the multiplication $x^2$
  on $L^2(\mathbb R^+, \mu)$
  is the positive half line $\mathbb R^+=[0, \infty)$. 
\end{prop+}

We shall need rescaled Wilson polynomials 
$\widetilde \omega_m(x^2): = \widetilde \omega_m((\frac{x} 3)^2)$. 
The factor of $3$ in our case is closely related to 
our imbedding 
of $SL(2, \mathbb R)$ in 
 $Sp(2, \mathbb R)$ via triple 
symmetric tensor product $S^3\mathbb R^2=\mathbb R^4$. 
We denote the corresponding measure by
$$
d\widetilde \mu(x) =d \mu_{a, b, c}(\frac x3).
$$

\begin{coro+}
  \label{c-d-h}
  Let $a\ge 0, b>0, c>0$, $A_m, C_m$ be as above. 
  The polynomials  $\widetilde \omega_m(x^2)  $
form an 
 orthogonal basis, 
 $L^2(\mathbb R^+, \widetilde\mu)$, 
\begin{equation}
    \label{w-m-sq-var}
    \Vert \widetilde
    \omega_m 
  \Vert^2 
  =\frac{m!(b+c)_m }{(a+b)_m(a+c)_m}, 
\end{equation}
  and satisfy the following recurrence relation 
  \begin{equation}
    \label{recur-Hahn}
    -\frac 1{3^2} (x^2 +d) 
\widetilde \omega_{m}(x^2)  = A_m  \widetilde \omega_{m+1}(x^2) 
  -(A_m +C_m -a^2 +\frac d{3^2})  \widetilde \omega_{m}(x^2) 
  +C_m   \widetilde \omega_{m-1}(x^2), 
  \end{equation}
  for any $d\in \mathbb R$. 
\end{coro+}
The introduction of $d$ is purely for practical convenience 
for identifying the spectrum of some abstract operator with 
the multiplication operator $x^2 +d$ in Section 3. 

\section{The  group 
  homomorphism $\iota: G_0=SU(1, 1)\to 
  G$ and the restriction of representations of $G$ to $G_0$}

We find the explicit realization of the principal $G_0=SU(1, 1)$
($=SL(2, \mathbb R)$) subgroup in the symplectic group $G$
and find explicit formulas for pull-back
of the metaplectic representation to $G_0$.

\subsection{The  homomorphism $\iota: G_0\to 
  G$}
The defining action of $SL(2, \mathbb C)$  on $\mathbb C^2$
preserves the complex symplectic form $dz_1\wedge dz_2$.
Its symmetric power representation $\iota:
SL(2, \mathbb C)\to SL(2n, \mathbb C)$ on
 $S^{2n-1}\mathbb C^2=\mathbb C^{2n}$
preserves   the symplectic form $S^{2n-1} (dz_1\wedge dz_2)$.
This defines a  group homomorphism
$\iota: SL(2, \mathbb C)\to Sp(n, \mathbb C)$.

It is immediate that $\iota(G_0)$,
$G_0=SU(1, 1)$,  preserves also  the Hermitian
form $\langle \cdot, \cdot\rangle_{n, n}=
S^{2n-1}\langle \cdot, \cdot\rangle_{1, 1}$ of
signature $(n, n)$. Thus  we have
the group homomorphism
\begin{equation}
 \label{def-iota}
\iota: G_0=SU(1, 1)\to U(n, n)\cap Sp(n, \mathbb C)=Sp(n, \mathbb R)=G.  
\end{equation}

We can find explicitly the image of $\iota$
of $SU(1, 1)$ and of the Lie algebra
$sl(2, \mathbb C)$. Let $m=2n-1$ and
    realize 
 $V=S^{m}\mathbb C^2$
    as  the space  homogeneous polynomials $f(x, y)$ degree of $m$ on 
    $\mathbb C^2$, 
    $g\in GL(2, \mathbb C): f(x, y)\mapsto f((x, y) g)$. 
    The space     $V$ is equipped with  the dual of Hermitian 
    form     $ S^{2n-1}\langle \cdot, \cdot\rangle_{1, 1}$ and 
    dual of the  symplectic form $  S^{2n-1} (dx\wedge dy)$, 
    namely $ \Omega(p, q)=(
    \partial_{x_1}
    \partial_{y_2}-
        \partial_{y_1}
        \partial_{x_2})^m (p(x_1, y_1)q(x_2, y_2)) $. 
   We fix an $\langle \cdot, \cdot\rangle_{n, n}$-orthonormal
    and $\Omega(\cdot, \cdot)$-symplectic basis of $V$,
    $$
    p_k=\sqrt{\binom m {2(k-1)}} x^{m-2(k-1)}{y^{2(k-1)}}
=\sqrt{\binom m {2(n-k)+1}} x^{m-2(k-1)}{y^{2(k-1)}},
$$
$$   
    \quad q_k=
\sqrt{\binom m {2(k-1)}}
        y^{m-2(k-1)}{x^{2(k-1)}}
=\sqrt{\binom m {2(n-k)+1}}
        y^{m-2(k-1)}{x^{2(k-1)}}, \quad 1\le k\le n;
       $$
namely they satisfy
$$
\langle    p_k, p_j\rangle_{n, n}
=-\langle    q_k, q_j\rangle_{n, n} =\delta_{kj},\quad
\Omega(p_k, p_j)=\Omega(q_k, q_j)=0, \quad
\Omega(p_k, q_j)=\delta_{kj}, \, 1\le k, j\le n.
$$

\begin{prop+}
  \label{sl-in-sp}
    Let $ k_\theta$ be the diagonal matrix
    $$   k_\theta 
    =\diag(
e^{(2n-1)i\theta}, e^{(2n-1-4)i\theta}, 
\cdots, e^{-(2n-3)i\theta}). 
$$
Then $\iota(u_\theta) =  k_\theta$ and 
the Lie algebra $\fg^{\mathbb C}$-elements
$$
H:=\iota(h),  E^{\pm }:=\iota(e^{\pm}), 
E=\iota(e) = E^+ + E^-
$$
are given by
$$
H=
\begin{bmatrix} D& 0\\ 0& -D 
\end{bmatrix},  
\quad E^+=\begin{bmatrix} 0& B\\ C& 0 
\end{bmatrix}, \quad
E^-=(E^+)^t, \quad E=
\begin{bmatrix} 0& B+C\\ B+C& 0 
\end{bmatrix}
$$
where
$$D=\diag(2n-1, 2n-5, \cdots, -(2n-3)),$$ 
$B$  is skew diagonal and symmetric
and $C$ lower skew diagonal and symmetric,
$$
B=\begin{bmatrix} 0& \dots & 0  & \beta_1
  \\
  0& \dots &   \beta_2 &0
  \\
  \vdots & \udots &  \vdots &  \vdots  
  \\
  \beta_n & \dots & 0& 0
\end{bmatrix}, \quad
C=\begin{bmatrix} 0& \dots & 0  & 0
  \\
  0& \dots &  0  &\gamma_2
  \\
  \vdots & \vdots &  \udots &  \vdots  
  \\
 0 & \gamma_n & 0 &  0
  \end{bmatrix}, 
  $$
  with $$
  \beta_k =\sqrt{(2k-1)(2(n-k)+1)}, \quad
  \quad 1\le k\le  n,$$
  $$
  \gamma_k=2 \sqrt{(k-1)(n-k+1)},
  \quad 2\le k\le  n; $$
  that is 
  $B=(b_{jk})$,  $C=(c_{jk})$ with $$
  b_{jk}=\beta_k\delta_{j, n-k+1}, \quad
  c_{jk}=\gamma_k\delta_{j, n-k+2}, \,
  1\le k, j\le n
  $$ in term
  of the Kronecker symbol $\delta_{pq}$.
\end{prop+}

\begin{proof}          The action of $k_\theta=\iota(u_{\theta})$ is  diagonal
       and is found  immediately, so is $H=\iota(h)$. Performing
       $\frac{d}{dt}|_{t=0}$ on the action
        $\iota(\exp(t e^+)) =
       \iota(\begin{bmatrix} 1 & t \\ 0 & 1
         \end{bmatrix}
         ):  f(x, y)\mapsto f(x, tx+y)$
         we find  $$E^+p_1=0,$$
         \begin{equation*}
           \begin{split}
                      E^+p_k&=\pi(e^+)p_k 
         =(2k-2) 
         \sqrt{\binom m {2(k-1)}}x^{m-2(k-1) +1}{y^{2(k-1)-1}}
         \\
         &=\frac{(2k-2) 
           \sqrt{\binom m {2(k-1)}}
         }
         {
\sqrt{\binom m {2(k-1)-1}}
} q_{n-(k-2)} \\
&
         =2\sqrt{(k-1)(n-k+1)}q_{n-(k-2)} 
           \end{split}
         \end{equation*}
         for $k\ge 2$,
         and
         $$
         E^+ q_k=
         \sqrt{(2k-1)(m-2k+2)} p_{n-k+1}, \quad 1\le k\le n.
         $$
         The element $E^-=\pi(\iota(e^-))=(E^+)^t$ by our choice
         of the basis.
\end{proof}

  We remark that the 
  homomorphism of $SL(2, \mathbb C)$
  into $Sp(n, \mathbb C)$ realizes also
 $SL(2, \mathbb R)$
as a subgroup in a different real form  $Sp(n, \mathbb R)$
in   $Sp(n, \mathbb C)$, and it is the principle $SL(2, \mathbb R)$-subgroup 
of $G=Sp(n, \mathbb R)$; see \cite{Kostant} for the general study 
of principle $SL(2, \mathbb R)$-subgroup in real split groups. 
This specific case
of $SL(2, \mathbb R)$
in $G=Sp(n, \mathbb R)$
is also of interests in
topology \cite{BIW}.

\subsection{The induced action of Lie algebra 
   of $\fsl(2)$ on $\mathcal F$ for  $n=2$. }

 We  shall find polynomials $p$
 that are annihilated by the differential operator
 $(B\partial,  \partial)$,  and we shall call them
 $(B\partial,  \partial)$-{\it harmonic}, i.e.
 \begin{equation}
   \label{B-hmc}
   B(\partial, \partial) p(z) =0.
 \end{equation}

 When $n=2$ the matrices $k_\theta, B, C, D$ are of the
 form
         \begin{equation}
           \label{bc-n=2}
    k_\theta 
    =\diag(e^{3i\theta}, e^{-i\theta}), \,
           B=\begin{bmatrix} 0&    \sqrt 3
  \\
 \sqrt 3&  0
\end{bmatrix}, \,
C=\begin{bmatrix} 0& 0   \\
  0& 2
\end{bmatrix}, \,
D=\diag(3, -1).
  \end{equation}
  This combined with  Lemma \ref{Lie-alg-mu}
  gives the explicit form for the
  action $\Lambda\circ \iota$ of
  of  $\fsl(2)$. With some abuse of notation
  we denote it also by   $\Lambda$, and we find
         \begin{equation}
           \label{hee-action}
                 \Lambda(H)= -1 -(3z_1\partial_1-z_2\partial_2),\,
  \Lambda(E^+)
  =-{\pi}^{-1} \sqrt{3} \partial_1\partial_2 +\pi z_2^2,\, 
    \Lambda(E^-)
    =-{\pi}^{-1}\partial_2^2 +\pi \sqrt{3} z_1 z_2.
  \end{equation}
  We can also prove independently that this defines
  a representation of $\fsl(2)$ on the Fock space
  without using the metaplectic representation. Indeed
  we have for any $f$, 
  \begin{equation*}
    \begin{split}
  [\Lambda(E^+),   \Lambda(E^-)]f
&  = [-{\pi}^{-1} \sqrt{3} \partial_1\partial_2
  +\pi z_2^2,
  -{\pi}^{-1}\partial_2^2
  +\pi \sqrt{3} z_1 z_2
  ]f
  \\
  &  =-3 [\partial_1\partial_2, z_1z_2] f
  -[z_2^2, \partial_2^2]f
  \\&=-3 \(z_1\partial_1 f +z_2\partial_2 f +f\)
  +2 f + 4 z_2 \partial_2 f 
    =-3z_1\partial_1 f +z_2 \partial_2 f  -f
  = \Lambda(H) f,
    \end{split}
  \end{equation*}
  \begin{equation*}
    \begin{split}
  [\Lambda(H),   \Lambda(E^+)]f
  &=  [-3z_1\partial_1
  + z_2\partial_2,
  -{\pi}^{-1} \sqrt{3} \partial_1\partial_2
  +\pi
  z_2^2
  ]f \\
&  = {\pi}^{-1} 3\sqrt{3}
  [z_1\partial_1, \partial_1\partial_2]f
  -{\pi}^{-1} \sqrt{3}
[z_2\partial_2,
    \partial_1\partial_2
    ]
    + \pi[z_2\partial_2, z_2^2]
    \\
    &= -{\pi}^{-1} 3\sqrt{3}\partial_1\partial_2f
    +{\pi}^{-1} \sqrt{3}\partial_1\partial_2f
    +2\pi z_2^2 f
= -{\pi}^{-1} 2\sqrt{3}\partial_1\partial_2f
+2\pi z_2^2 f
\\
&=2(-{\pi}^{-1} \sqrt{3}\partial_1\partial_2f
+\pi z_2^2 )=2\Lambda(E^+)f,  
    \end{split}
  \end{equation*}
  and 
  \begin{equation*}
    \begin{split}
  [\Lambda(H),   \Lambda(E^-)]f
&  = [-3z_1\partial_1
  + z_2\partial_2, -\pi^{-1}\partial_2^2 + \pi\sqrt{3}z_1z_2]f
  \\
  &=-3\pi\sqrt{3}
  [z_1\partial_1, z_1z_2]f
  -\pi^{-1}[z_2\partial_2, \partial_2^2]f
  +\pi\sqrt{3}[ z_2\partial_2, z_1z_2]f
  \\
  &=-3\pi\sqrt{3}z_1z_2 f + 2\pi^{-1}\partial_2^2 f
  +\pi\sqrt{3} z_1z_2 f
  =-2\Lambda(E^-) f.
    \end{split}
  \end{equation*}

  \begin{rema+}  It was pointed out to me by the anonymous referee
    that $\fsp(2, \mathbb R)
    $ has two principal
    nilpotent orbits under the adjoint action of
    $Sp(2, \mathbb R)$.
    In terms of the simple roots
    $ \{\epsilon_1-\epsilon, 2\epsilon_2\}$
    for the Lie algebra     $\fsp(2, \mathbb R)$, the
    complex principal nilpotent orbit
    is determined, up to scalars,  by the element $H=2\epsilon_1 +
    (\epsilon_1+\epsilon_2) =3 
    \epsilon_1 +\epsilon_2$   in a standard $\fsl(2)$-triple $\{H,
    E^{\pm}\}$,     
    namely the diagonal matrix $\diag(3, 1,
    -3, 1)$. The nilpotent element $E^+$ has two
    possible forms,
    $$
    E^+: q_1\to p_2\to q_2\to p_1\to 0, \quad
        E^+: p_1\to q_2 \to p_2\to  q_1\to 0.
    $$
    Our principal $\fsl(2, \mathbb R)$ subalgebra is conjugated
    to the first one; the two nilpotent matrices are conjugated
    by a matrix exchanging the symplectic forms $\Omega$ and $-\Omega$.
\end{rema+}
    


  \section{The complete
    decomposition
    for $( \Lambda, Mp(2, \mathbb R)    )$
  under $SU(1,1)$}

\subsection{Orthogonal basis for $U(1)$-weight vectors}
Let $n=2$.
We shall find orthogonal polynomials
in $\mathcal F=\mathcal F(\mathbb  
C^2)$ of a fixed $U(1)=\iota(U(1))$-weight.
Let $\mathcal P\subset \mathcal F$ be the polynomial space.
All weights refer to the metaplectic representations
under $U(1)$
unless otherwise explicitly stated.

Denote
\begin{equation} \label{invariant-I}
  I(z)=\frac {\pi^2}{3\sqrt 3}z_1 z_2^3.\end{equation}
(The coefficient
    is chosen to simplify the expression for the solutions of the
    equation (\ref{hwv-eq}) below.)
    Then $I$ generates all invariants
      in the polynomial space $\mathcal P$
   of the defining action of $\iota(U(1))$.  
Recall $(\alpha)_m
    =\alpha (\alpha +1)\cdots (\alpha +m-1)$,  the Pochhammer symbol.

    The
 orthogonal  basis vectors $  \{z_1^{m_1} z_2^{m_2}\}$ of $\mathcal F $
 are weight vectors of $U(1)$
    of weights $-1 -3m_1 +m_2$, and modulo $3$
    they are of the form $\mu=-3k-1, -3k, -3k+1$. We denote
    $\mathcal F|_{U(1)}^{\mu}$ the subspace of all
    weight vectors of weight $\mu$.
    Then if $k\ge 0$ we have
    \begin{equation}
      \label{w-v-1}
    \mathcal F|_{U(1)}^{-3k-1} =\sum^{\oplus}_{l\ge 0}\mathbb C(I^l(z) z_1^k),      
    \end{equation}
    \begin{equation}
      \label{w-v-2}   
    \mathcal F|_{U(1)}^{-3k} =\sum^{\oplus}_{l\ge 0}\mathbb C(I^l(z) z_1^kz_2),
    \end{equation}
    and
    \begin{equation}
      \label{w-v-3}
      \mathcal F|_{U(1)}^{-3k+1} =\sum^{\oplus}_{l\ge 0}\mathbb C(I^l(z) z_1^kz_2^2);
    \end{equation}
if $k<0$ we  then replace $z_1^k$ by $z_2^{-3k}$ in the above formulas. 

We compute the norm of weight vectors.
\begin{lemm+}
  \label{norm-1}
  The square norms      $ \Vert I^mz_1^k\Vert^2_{\mathcal F}$
and        $ \Vert I^mz_2^k\Vert^2_{\mathcal F}$
     are given by
      $$
      \Vert I^m z_1^k\Vert^2=\frac{1}{\pi^k}
     (m!)^2 (m+1)_k (\frac 23)_m(\frac 13)_m,
      $$
            $$
      \Vert I^m z_2^k\Vert^2=\frac{1}{\pi^k}
      (m!)^2 (3m+1)_k (\frac 23)_m(\frac 13)_m. 
      $$
  \end{lemm+}

  \begin{proof} Using (\ref{monom-norm}) we see
    that
    $$
    \Vert I^m z_1^k\Vert^2= (\frac {\pi^2}{3\sqrt 3})^{2m}
    \Vert z_1^{m+k}  {z_2^{3m}}
    \Vert^2
    =\frac {\pi^{4m}}
    {3^{3m}}
    \frac 1{\pi^{4m+k}}    (m+k)! (3m)!
    =\frac {\pi^{4m}}
    {3^{3m}}
    \frac 1{\pi^{4m+k}}   m! (m+1)_k (3m)!.
    $$
    Now
    $$(3m)! = (3m) (3m-1) (3m-2)
    \cdots
    3\cdot 2\cdot 1
    = 3^{3m} m (m-\frac 13)(m-\frac 23)
    \cdots    1\cdot \frac 23 \cdot \frac 13
    = 3^{3m} m! (\frac 23)_m(\frac 13)_m,
    $$
    and this proves the first formula.

    Similarly
    $ \Vert I^mz_2^k\Vert^2_{\mathcal F}$
    is
    $$
\frac {\pi^{4m}}
    {3^{3m}}
    \frac 1{\pi^{4m+k}}
    m! (3m+k)!
    =\frac 1{\pi^{k}}
    \frac 1 {3^{3m}} m! (3m+1)_k (3m)!
$$
with  $(3m)!$
being computed as above.
\end{proof}

      \subsection{Irreducible decomposition of $(\Lambda, Sp(n,
        \mathbb R))$}
Our main result is the following
  
\begin{theo+}
    \label{main}
    The  decomposition of      $(\Lambda, \mathcal F, Mp(2, R)
    )$ under  $SU(1, 1)$ is given by
    \begin{equation}
      \label{eq:1}
(\Lambda,  Mp(2, R))|_{SU(1, 1)}
\cong \int_0^\infty \sigma_{i\lambda, +} d\tilde\mu_{\frac 12, \frac
  16, \frac 56}(\lambda)
\oplus 
\int_0^\infty \sigma_{i\lambda, -} d\tilde\mu_{0, \frac 13, \frac 23}(\lambda)
\oplus 
\bigoplus_{k=1}^{\infty} \sigma_{-3k-1},
\end{equation}
were $\tilde\mu_{a, b, c}$ is the orthogonality measure for the 
continuous dual Hahn polynomials in Corollary \ref{c-d-h}. 
  \end{theo+}
  \begin{rema+}    The measure 
    $\tilde\mu_{a, b, c}$ can be written as $|c_l(\lambda)|^{-2}b(\lambda)^{-1}$
    where $c_l(\lambda)$ is the Harish-Chandra $c$-function
    for  line bundles with parameter $l$ \cite{Sch}
    and $b(\lambda)$ is the symbol of a Berezin transform \cite{gz-bere-rbsd}.
    We can follow the method of Berezin transform
    to study the decomposition above, however it requires
    different realization of the metaplectic representation
    and it is less effective in finding the discrete components;
    see Remark
    \ref{rema-be}.
  \end{rema+}
    
We shall find the spectral
    decomposition of the Casimir element
     $C    $      on $\mathcal F$.
     We have
 $\iota(C)=E^+E^-
 +E^- E^+ +\frac{H^2}2$ as element
 in the enveloping algebra of $\fsp(n, \mathbb C)$, 
 and 
 $$
 \Lambda (\iota(C))=\Lambda(E^+)
 \Lambda
 ( E^-)
 +\Lambda (E^-)
 \Lambda( E^+) +\frac{\Lambda(H)^2}2
 $$ as
 operator on $\mathcal F$. To ease
 notation we write it just as $C$.
The decomposition is  done in several steps. 

By general abstract theory the operator $C$
has a well-defined self-adjoint extension on   $\mathcal F$
and on any subspace
$\mathcal F|_{U(1)}^{k}$ of fixed weight $k$.
 \begin{lemm+} 
    \label{3.3}
   The Casimir element $-C$ on 
   $\mathcal F^{-1}|_{U(1)}$
   has continuous spectrum $[\frac{1}{2}, \infty)$.
 \end{lemm+}

 \begin{proof}   The Casimir operator    $C=E^+E^-
   +E^- E^+ +\frac{H^2}2 = 2E^+E^- - H +
   \frac{H^2}2$ acts on any element $f\in \mathcal P|_{U(1)}^{-1}$
   of weight $-1$, $\Lambda(H)f=-f$,  as
    $$
    Cf= 2 \Lambda( E^+) 
    \Lambda (E^-) f+ f +\frac 12 f =
     2 \Lambda( E^+) 
    \Lambda (E^-)f + \frac 32 f. 
    $$
The subspace 
   $ \mathcal F|_{U(1)}^{-1}  $
   has an orthogonal basis
   $W_m(z)= \frac{I^m (z)}
{(\frac 13)_m (\frac 23)_m 
   }$
   with 
   \begin{equation}
      \label{W-m-sq}     
   \Vert   W_m\Vert^2 =\frac{1}{\pi} \frac{(m!)^2}{
     (\frac 13)_m (\frac 23)_m 
   } 
   \end{equation}
   by Lemma   \ref{norm-1}.
   We compute the action of $\Lambda(E^+)
   \Lambda(E^-)$
   on  $\{W_m\}$. We have,
$    \pi  \Lambda 
   ( E^-)=-\partial_2^2 +\pi^2\sqrt 3 z_1 z_2$, and
   by straightforward computations,
   \begin{equation}\label{lambda-e-I-m}
     \begin{split}
        \pi  \Lambda 
   ( E^-) 
   I_m(z) 
&=  (-\partial_2^2 +\pi^2 \sqrt{3} z_1 z_2) 
\left(\frac{\pi^2}{3\sqrt 3} z_1z_2^3 
\right)^m \\
&=-(3m)(3m-1) I^m(z) z_2^{-2}
+ 3^2I^{m+1}(z) z_2^{-2}\\
&=-3^2m(m-\frac 13)I^m z_2^{-2}
+  3^2I^{m+1} (z)z_2^{-2}.
   \end{split}            
   \end{equation}
Acting by $
 \pi  \Lambda 
( E^+) = -\sqrt 3 \partial_1\partial_2 +\pi^2 z_2^2$ we find,
writing $I(z)^{m} = I(z) I^{m-1}(z)
=\frac{\pi^2 z_1z_2^3}{3\sqrt 3} 
I^{m-1}
$, that
   \begin{equation*}
     \begin{split}
&\quad   \pi^2    \Lambda( E^+)  \Lambda( E^-) 
  I_m\\
  &  =   -\sqrt 3\partial_1\partial_2
\left( -3^2m(m-\frac 13)I^m z_2^{-2}
  +  3^2I^{m+1} (z)z_2^{-2}\right)
+\pi^2 z_2^2
\left( -3^2m(m-\frac 13)I^m z_2^{-2}
  +  3^2I^{m+1} (z)z_2^{-2}\right)\\
&=
  -\sqrt 3
  \left( -3^3 m(m-\frac 13)m(m-\frac 23)
    \frac{\pi^2}{3\sqrt 3}
    I^{m-1}
    + 3^3 (m+1)(m+\frac 13) \frac{\pi^2}{3\sqrt 3} I^m
  \right)\\
  &\quad   +\pi^2 3^2\left(
    -m (m-\frac 13) I^m +I^{m+1}\right)
  \\
  &=     3^2\pi^2  (\alpha_{m} I^{m+1} +\beta_{m} I^{m} +\gamma_m I^{m-1})
       \end{split}
   \end{equation*}
   with the leading coefficient
   $   \alpha_m=1,
   $
   the last coefficient
   $$\gamma_m= m^2(m-\frac 13) 
   (m-\frac 23)
   $$
   and the middle
   $$
   \beta_m=-(m+1)(m+\frac 13) -m(m-\frac 13)=-(2m +m +\frac 13).
   $$
         Thus
         $$
 \frac 1{3^2}   \Lambda( E^+)
      \Lambda( E^-)
      I^m
      = I^{m+1}
            -(2m^2 + m +\frac 13   ) I^m
      +m^2 (m-\frac 13) 
      (m-\frac 23)      I^{m-1}.
      $$
      Writing in terms
      of $ W_m(z) =\frac{I^m (z)}{
           (\frac 13)_m (\frac 23)_m }
         $ and  using  $  (\frac 13)_{k+1} (\frac 23)_{k+1}
         = (\frac 13)_{k} (\frac 23)_{k} (k+\frac{1}3) ( k+\frac 23)$ for $k=m, m+1$,          
         this becomes
      $$  \frac 1{3^2}
      \pi( E^+)
      \pi( E^-)
      W_m
      =(m+\frac 13)(m+\frac 23)W_{m+1}
      -(2m^2 + m +\frac 13   )  W_{m}
    +m^2      W_{m-1}.
    $$
This is exact the same recursion relation
as    (\ref{recur-Hahn})
    for the continuous dual Hahn polynomials
    for
    $$a=0,\,  b=\frac 13, \, c=\frac 23, \, d=1;
    $$ 
   moreover  square norms (\ref{W-m-sq}) and
   (\ref{w-m-sq-var}) are same. In other words,
   $W_m \to \widetilde \omega_m$ is a unitary operator
   intertwining $-\Lambda( E^+)    \Lambda( E^-)$
with the multiplication operator
by $1+x^2$  on the space $L^2(\mathbb R^+, \widetilde \mu)$.
Thus $-\Lambda( E^+)    \Lambda( E^-)$
    has continuous spectrum $[1, \infty)$. But $
        -Cf= -2 \Lambda( E^+) 
        \Lambda (E^-)f - \frac 32 f$ so $-C$ has
        spectrum $[\frac 12, \infty)$. 
    This finishes the proof.
  \end{proof}

  Note that the recursion formula (\ref{lambda-e-I-m}) can also
  be obtained by just finding the leading term $3^2 I^{m+1}z_2^{-2}$
  and
  by using the unitarity of $\Lambda$.

  \begin{lemm+}
    \label{3.4}
   The Casimir element $-C$ on 
   $\mathcal F^0|_{U(1)}$
   has continuous spectrum $[\frac{1}{2}, \infty)$.
 \end{lemm+}   
 \begin{proof} This is proved by similar computations
   as the above lemma.
   We consider the action of $C$ on $\mathcal P^0|_{U(1)}\subset 
    \mathcal F^0|_{U(1)}  $, 
    the subspace of polynomials $f$ of weight $0$.
    If $f\in \mathcal P|_{U(1)}^0$ then      $\Lambda(H) f=0 $
and     $$
Cf= 2 \Lambda(E^+)
\Lambda(E^-)f.
    $$
   The  space
   $ \mathcal F|_{U(1)}^0  $ has an orthogonal
   basis given by  the polynomials
   $W_m(z)= \frac{I^m(z) z_2}
{(\frac 23)_m (\frac 43)_m 
   }$,   with square norms  computed in  Lemma   \ref{norm-1},
   $$
   \Vert   W_m
   \Vert^2 =
   \frac{1}{\pi} \frac{(m!)^2}{
     (\frac 23)_m (\frac 43)_m 
   }. 
   $$
   We compute the action of $C$ on $I^m z_2 $
   and find
$$
\frac 1{3^2} \Lambda(E^+)
\Lambda(E^-) (I^m z_2)
=I^{m+1}(z) z_2 - (2m^2 +2m +\frac 23)  (I^m z_2)
+ m^2 (m-\frac 13)(m +\frac 13)  (I^{m-1} z_2).
$$
Written in terms of   $ W_m  $
this is
$$
\frac 1{3^2} \Lambda(E^+)
\Lambda(E^-) W_m
= (\frac 23 +m)
(\frac 43 +m)  W_{m+1}
-(2m^2 +2m +\frac 23) W_{m}
+m^2  W_{m-1}.
$$
This is the same recursion formula for continuous Hahn polynomials in
for
$$
a=\frac 12,\,  b=\frac 16,  \, c=\frac 56, \, d=\frac 14,
$$
and it follows that the
$\Lambda(E^+)\Lambda(E^-)$ with the multiplication operator
by $-(d +x^2) =-(\frac 14 + x^2) $ on
$L^2(\mathbb R^+, \widetilde \mu)$.
This gives the spectrum of
$\Lambda(E^+)\Lambda(E^-)$ and the Casimir element as claimed.
\end{proof}

We study then the appearance
of lowest or highest weight representations.

\begin{lemm+}
    \label{exist-l-h}
    There exists no   lowest weight unitary representation of $SU(1, 1)$
    in   $(\Lambda, \mathcal F)$.
    The only highest weight unitary representations of $SU(1, 1)$
    in  $(\Lambda, \mathcal F)$ are of the form $\sigma_{-3k-1}$,
    $k>0$, and they all appear of multiplicity one.
 \end{lemm+}   

 \begin{proof}
   Let $f$ be a lowest weight vector, i.e., $f$
    is a weight vector and $\Lambda(E^-)f=0$. Now
$$
E^-=(E^+)^t=\begin{bmatrix} 0& C\\ B& 0 
\end{bmatrix}
$$
with $B $ and $C$ as in (\ref{bc-n=2}).
Thus the lowest weight vector condition for $f$  becomes 
$$
-\frac1{\pi} (\partial_2^2 f)(z)
 +\pi \sqrt 3 z_1 z_2f =0.
$$
Then the lowest
degree term $f_0$ of $f$ satisfies $\partial_2f_o=0$
and must be of the form $f_0=z_1^k$ or
$f_0=z_1^kz_2$ up to non-zero constants. But these vectors
have negative weights under $\Lambda(H)$
and any unitary lowest weight irreducible representation of $\fsu(1, 1)$ has positive
lowest weight, consequently lowest weight representations
do not appear.

Next we shall find all  highest weight vectors in the
    space $\mathcal F$. 
 So let $\sigma_\nu$ be an irreducible
    representation
    for the Lie algebra $su(1, 1)$ appearing in $\Lambda$ with 
    highest weight vector $f$ of weight $-\nu$, $\nu >0$.

Expand $f$    as sum of  homogeneous polynomials
    $$
f=\sum_{m=0}^\infty p_{m+k},$$
with $p_{m+k}$  of homogeneity $m+k$ and $k$ being the lowest degree.
The highest weight condition $\Lambda(E^+) f=0$ becomes
\begin{equation}
  \label{hwv-eq}
\frac 1{2\pi} (B\partial, \partial) f=
\frac {\pi}2 (Cz, z) f,
\end{equation}
i.e.
$$
\frac 1{2\pi} (B\partial, \partial) p_{m+k}=
\frac {\pi}2 (Cz, z) p_{m+k-4}, m=0, 1, \dots.
$$
Thus $p_{m+k-4}$ determines $p_{m+k}$ up to $B$-harmonic polynomials
of  degree $m+k$.
Hence the sum $\sum_{l=0}^\infty p_{k+4l}$ in $f$
satisfies the equation (\ref{hwv-eq}) and it starts with
the lowest degree $p_k$ of $f$. But since we
are finding all solutions we can assume $f=\sum_{l=0}^\infty p_{k+4l}$.
To ease notation
we write this as $f=\sum_{l=0}^\infty f_{l}$, $f_l=p_{k+4l}$.

We have then
\begin{equation}\label{recur-f}
\frac{\sqrt 3}{\pi^2} \partial_1\partial_2 f_{l+1}=
 z_2^2 f_{l}, l=0, 1, \cdots, \end{equation}
and  the leading term $f_0$ satisfies
$\partial_1 \partial_2f_0=0$, namely
    $2\sqrt 3\partial_1\partial_2 f_0=0$. So $f_0= z_1^k$,  $k\ge 0$
    or $f_0=z_2^k, k>0$.

    If $f_0=z_2^k$ then
    $\Lambda(H)f_0= (-1+k) f_0$,
    it is of nonnegative weight $-1+k$ and this contradicts our
     assumption that $f$ is of negative weight, so this case is excluded.

     Let   $f_0=z_1^k$, $k\ge 0$. The weight of $f$ is $\nu=-1-3k$.
     We claim that $f$ is uniquely solved by
    \begin{equation}\label{solu-f}
    f(z)=k!z_1^{k}\sum_{l=0}^\infty \frac{I(z)^{l}}{l! (l+k)!}.
    \end{equation}

    Now the action
    of $U(1):=\iota(U(1))$ on the space of
    $\mathcal P_m$ of homogeneous polynomials
    of degree $m$ is multiplicity free for any $m$. Thus that
    $f_l$ is of weight $-k-1$
    and of degree $k+4l$
    implies that
    $$
    f_l
    =a_l I(z)^l z_1^k.
    $$
    Note that
    \begin{equation*}
      \begin{split}
\frac{\sqrt 3}{\pi^2} \partial_1\partial_2 
(I(z)^{l+1} z_1^k)
&=\frac{\sqrt 3}{\pi^2} 
(\frac{\pi^2}{3\sqrt 3})^{l+1}
\partial_1\partial_2 ( z_1^{l+k+1}z_2^{3(l+1)})
  \\
  &
=\frac{\sqrt 3}{\pi^2} 
(\frac{\pi^2}{3\sqrt 3})^{l+1} (l+k+1) 3(l+1)
( z_1^{l+k}z_2^{3l +2})
  \\
  &
  =  (\frac{\pi^2}{3\sqrt 3} z_1 z_2^3)^l (l+k+1)(l+1) z_2^2
  \\
  &= (l+k+1)(l+1) z_2^2    I(z)^l.
      \end{split}
    \end{equation*}
    The recursive relation  (\ref{recur-f}) becomes
    $$
    a_{l+1} =\frac{1}{(l+k+1)(l+1)}a_l
    $$
    and     is solved by
    $$
    a_l=\frac{k!}{l! (k+l)!}
    $$
This proves our claim.   

    Finally we prove that $f\in \mathcal F$ if $k>0$.
    The square norm of $f$ is
    \begin{equation*}
      \begin{split}
(  k! )^2 \sum_{l}\frac 1{(l!^2 ((k+l)!)^2} \Vert  I^l z_1^k\Vert^2
   & =    \sum_{l}\frac 1{l!^2 ((k+l)!)^2}
   (l!) (l+k)! (\frac 23)_l(\frac 13)_l
   \\
   &
=  \sum_{l}\frac {(l+k)! (\frac 23)_l(\frac 13)_l}
{l!((k+l)!)^2}
   \\
   &
   = k!\sum_{l}\frac {(\frac 23)_l(\frac 13)_l}
   {(1)_l(k+l)!}
   =
   k!
   \frac
      {
    \Gamma(k+1) 
   }
{\Gamma(\frac 23)\Gamma(\frac 13)}
   \sum_{l}\frac {\Gamma(\frac 23+l)\Gamma(\frac 13+l)}
   {
     \Gamma(l+1)\Gamma(k+l+1) 
   }   \\
   &
   \asymp \sum_l \frac 1{l^{k+1}}        
      \end{split}
    \end{equation*}
which is convergent iff $k\ge 1$. Here we have
used the known fact that
$$
\frac {\Gamma(a+l)}
{\Gamma(b+l)}\asymp \frac 1{l^{b-a}}.
$$
This finished the proof.\end{proof}
  
We can now prove Theorem   \ref{main}. 
\begin{proof}
It follows from  Lemmas  \ref{3.3}, \ref{3.4}
and \ref{exist-l-h} that right hand side
in     (\ref{eq:1}) are subrepresentation of
$(\Lambda,  \mathcal F, SU(1, 1))$. We
have to prove that it gives full decomposition.
The proof relies on some general and abstract arguments.
Let $\mathcal D\subset \mathcal F$ the
the subspace of $SU(1, 1)$-smooth vectors in
$\mathcal F$ and $\mathcal D^{\mu}\subset \mathcal F^{\mu}$ be
the subspace of all vectors of weight $\mu=-3k-1, -3k, -3k+1$.
The decomposition of $\mathcal F$ is obtained
by the spectral decomposition of the Casimir operator $C$
on $\mathcal D^\mu$ and its self-adjoint extension.
The operator $E^{-} $  has $\Ker E^{-}|_{\mathcal D}=0$, and all
elements in $\Ker E^{+}|_{\mathcal D}$ are found,  in Lemma
\ref{exist-l-h}. Moreover
the space
$\Ker E^{+}|_{\mathcal D}$ generates the discrete series
representation $\sigma_{-3k-1}$ in $\mathcal F$.
Let $\mu=-6k-1$, $k>0$. The operator
$(E^+)^{3k}$ defines an injective operator
$\mathcal D^{\mu }\cap (\Ker (E^+)^{3k})^{\perp}\to
\mathcal D^{-1} \subset \mathcal F^{-1}$ and it intertwines
the Casimir operator. The Casimir operator $C$
on  $\mathcal F^{-1}$ has only continuous spectrum
and thus  it has only continuous spectrum on
$\mathcal F^{\mu }\cap (\Ker (E^+)^{3k})^{\perp}$. Similar considerations
apply to the subspaces
$\mathcal F^{\mu }$ of all weights. This
give the spectral decomposition of $C$ on
$\mathcal F^{\mu}$ and then
on $\mathcal F$, the spectral measure being given in
Lemma \ref{3.3} and \ref{3.4}. This completes the proof.
\end{proof}

The  even part 
 $(Mp(2, R), \Lambda^+, \mathcal F_+)$
of
the metaplectic representation
consists of even holomorphic functions,
the corresponding discrete components
correspond to the
the highest weight vector (\ref{solu-f}) with $k=2l$ even.
We have thus

\begin{coro+}
  The discrete components appearing in 
  the restriction of
the even part of the metaplectic representation
$(\Lambda^+, \mathcal F_+, Mp(2, R))$
  to   $SU(1, 1)$ are precisely
  the highest weight representations
  of weights $-6l-1$, $l\ge 1$.
  \end{coro+}

  \begin{rema+}
    \label{rema-be}
The continuous part of the decomposition of the even part (as
holomorphic functions)
of the metaplectic
representation $\Lambda^+$ can 
also be found using its realization
as reproducing kernel space. We give a very
brief description.
The representation $\Lambda^+$ can 
be realize \cite{KV-invent} on the Siegel bounded symmetric
domain $D=G/K=Sp(2, \mathbb R)/U(2)$
as a Hilbert space $\mathcal H$ of holomorphic function with
reproducing kernel $\det (1-z\bar w)^{-\frac 12}$.
It is a subspace of holomorphic sections  
of a line bundle over the Siegel domain  
 $Sp(2, \mathbb R)/U(2)$.  
 Consider the
the hyperbolic disc 
$\Delta=SU(1, 1)/U(1)=\{z\in \mathbb C; |z|<1\} $ and 
the  Hilbert space $L^2(D, (1-|z|^2)^{-1} dm(z))$ 
with the natural $G_0=SU(1, 1)$ action as  
$$g\in G_0: f(z)\to  f((az+b)(cz+d)) (cz+d)^{-1}.$$
This is also the $L^2$-space of the half-bundle $\mathbb C (dz)^{\frac 12}$
of the cotangent bundle.
The pull-back of $\iota$,  (\ref{def-iota}), 
$$
R=\iota^\ast: \mathcal D\subset \mathcal H \to
L^2(D, (1-|z|^2)^{-1} dm(z))
$$ defines an intertwining map $R$
from a dense subspace of $H$ to 
$L^2(\Delta, (1-|z|^2)^{-1} dm(z))$. The operator $RR^\ast$
is a Berezin transform on $\Delta$ and its spectral symbol
can be computed using the method in \cite{gz-bere-rbsd}.
However the restriction map $R$ here is not 
injective and it will produce only the continuous part, and other
discrete
components can be
constructed using covariant differentiation; this
requires rather subtle and detailed computations.

  It is interesting to notice that  
the solution    (\ref{solu-f})  
for $k=0$  still defines a holomorphic function, which  
is a highest weight vector but not in the Fock space. This  
can also be explained by another fundamental fact about  
the Hardy space of holomorphic functions on the unit disc.  
This space
$(L^2(\Delta, (1-|z|^2)^{-1} dm(z)), G_0)$ has no discrete component. However the Hardy  
on the unit disc is   invariant  
under the same action but it is not a subspace 
$L^2(D, (1-|z|^2)^{-1} dm(z)).$  
So the solution  (\ref{solu-f}) plays the role of a Hardy space  
function whereas the restriction of the reproducing kernel
space to $D$ plays the role  of  $L^2(D, (1-|z|^2)^{-1} dm(z)).$

\end{rema+}

We observe that the proof of Lemmas
    \ref{3.3}  and   \ref{3.4}
actually also provides generalized eigenfunctions
for the Casimir operator as a sum of bi-orthogonal
polynomials
$$\Psi_{x}(z)=
\sum_{m=0}^\infty e_m (z) \widetilde e_m (x^2) 
$$
where $e_m(z)=\frac{W_m(z)}{\Vert W_m\Vert}$,
 $\widetilde e_m=\frac{\widetilde \omega_m(x^2)}{\Vert \widetilde
   \omega_m\Vert}$.
 It is easy to prove that this series
 is convergent point-wise for all $x\in \mathbb R$
 and $z\in \mathbb C^2$. We have
 not been able to compute the sum. There exist several
 formulas \cite[pp.~196-199]{Koe-Swa}
 for  generating functions
  $\sum_m c_m \widetilde \omega_m(x^2) t^m$
  of the polynomials $c_m\widetilde \omega_m(x^2)$
  however they are not related to the eigenfunctions
  $\Psi_{x}(z)$.


It might be possible to find all discrete
components for  $(\Lambda, \mathcal F, Mp(n, \mathbb R))$
for general $n\ge 3$ by refining
the techniques in Lemma \ref{exist-l-h}. It would also be an interesting problem to
  study the branching problems
  of minimal representations
  for other split real simple
  Lie groups, such as quaternionic representations
  of $G_2$, under its principal
  $SL(2, \mathbb R)$-subgroup.

\def\cprime{$'$} \newcommand{\noopsort}[1]{} \newcommand{\printfirst}[2]{#1}
  \newcommand{\singleletter}[1]{#1} \newcommand{\switchargs}[2]{#2#1}
  \def\cprime{$'$} \def\cprime{$'$} \def\cprime{$'$} \def\cprime{$'$}
\providecommand{\bysame}{\leavevmode\hbox to3em{\hrulefill}\thinspace}
\providecommand{\MR}{\relax\ifhmode\unskip\space\fi MR }
\providecommand{\MRhref}[2]{%
  \href{http://www.ams.org/mathscinet-getitem?mr=#1}{#2}
}
\providecommand{\href}[2]{#2}




\end{document}